\documentclass[12pt]{article}
\usepackage{amsmath}
\usepackage{verbatim}
\usepackage{amssymb}
\usepackage{tikz}
\usetikzlibrary{arrows,decorations.markings}
\newtheorem{thm}{Theorem}
\newtheorem{lemma}[thm]{Lemma}
\newtheorem{cor}[thm]{Corollary}
\newtheorem{prop}[thm]{Proposition}

\newcommand\e{\varepsilon}
\newcommand{\g}{\boldsymbol{\gamma}}
\newcommand{\al}{\boldsymbol{\alpha}}
\newcommand{\be}{\boldsymbol{\beta}}
\newcommand{\muu}{\boldsymbol{\mu}}
\newcommand{\om}{\boldsymbol{\omega}}
\newcommand{\et}{\boldsymbol{\eta}}
\newcommand{\si}{\boldsymbol{\sigma}}
\newcommand\Ss{\mathcal{S}^-\!}

\begin{document}
\title{On  $\infty$--Ground States in the Plane}
\author{Erik Lindgren, Peter Lindqvist}
\date{\empty}
\maketitle

\maketitle

\medskip
{\small \textsc{Abstract:} \textsf{We study $\infty$-Ground states in convex domains in the plane. In a polygon, the points where an $\infty$-Ground state does not  satisfy the $\infty$-Laplace Equation are characterized: they are restricted to  lie on specific curves, which are  acting as attracting  (fictitious) streamlines.  The gradient is continuous outside these curves and no streamlines can meet there. }}

\bigskip

\noindent {\small \textsf{AMS Classification 2000}: 35K65, 35P30, 35J70}


\section{Introduction}

The $\infty$--Ground state was defined in \cite{JLM} as a viscosity solution $u\in W^{1,\infty}_0(\Omega),\, u > 0$ in $\Omega$, of the equation
\begin{equation}\label{Ground}
  \mathrm{max}\Bigl\{\Lambda - \frac{|\nabla u|}{u},\,\sum_{i,j}\frac{\partial^{\,2} u}{\partial x_i\,\partial x_j}\frac{\partial u}{\partial x_i}\frac{\partial u}{\partial x_j}\Bigr\}\,=\,0.
\end{equation}
Here $\Omega$ is a domain in $\mathbb{R}^N.$ In fact, $u \in C(\overline{\Omega})$ and $u = 0$ on the boundary $\partial \Omega$. The equation is obtained as the limit of the Euler-Lagrange equations
\begin{equation}\label{peqn}
  \nabla\!\cdot\!\bigl(|\nabla u|^{p-2}\nabla u\bigr)\,+\,\lambda_p|u|^{p-2}u\,=\,0
  \end{equation}
  in the problem of minimizing the Rayleigh quotients
  \begin{equation*}
    \lambda_p\,\,=\,\,\underset{u\in W^{1,p}_0(\Omega)}{\mathrm{min}}\frac{\displaystyle\int_{\Omega}\!|\nabla u|^p\,dx}{\displaystyle\int_{\Omega}\!|u|^p\,dx}.
  \end{equation*}
  The positive minimizers are called $p$--Ground states. The operator  $\Delta_{\infty}u\equiv\sum\frac{\partial^{\,2} u}{\partial x_i\,\partial x_j}\frac{\partial u}{\partial x_i}\frac{\partial u}{\partial x_j}$ is the $\infty$-Laplacian, which was introduced by G. Aronsson in \cite{A1}. It is the formal limit of the $p$-Laplace operators $\Delta_pu\,\equiv\,\nabla\!\cdot\! (|\nabla u|^{p-2}\nabla u).$
  
  Some features which are typical for eigenvalue problems are preserved at the limit as $p \to \infty$. A most remarkable property is that the $\infty$--Ground state exists if and only if $\Lambda$ has the specific value
  \begin{equation*}
    \Lambda_{\infty}\,=\,\lim_{p \to \infty}\sqrt[p]{\lambda_p}\,=\,\frac{1}{\underset{x \in \Omega}{\mathrm{max}}\,\mathrm{dist}(x,\partial \Omega)}.
  \end{equation*}
  In other words, the $\infty$-eigenvalue $\Lambda_{\infty}$ is equal to the reciprocal value of the radius of the largest ball that can be inscribed in $\Omega$. Such a simple rule is not known even for $p =2$, when equation (\ref{peqn}) reduces to the Helmholz Equation $\Delta u + \lambda\,u\,=\,0.$ The proof in \cite{JLM} was based on the equation
  for $v = \log(u)$:
  \[
    \mathrm{max}\Bigl\{\Lambda -|\nabla v|,\,\Delta_{\infty} v + |\nabla v|^4\Bigr\}\,=\,0.
 \]
 
 Those solutions of (\ref{Ground}) that come as limits of the eigenfunctions $u_p$ in (\ref{peqn}) are called \emph{variational} $\infty$--Ground states.  There is always at least one obtained  in this way. It was shown by R. Hynd, Ch. Smart, and Y. Yu that the (normalized)  $\infty$--Ground state is not unique. Their counterexample in \cite{HSY} is a dumbbell shaped domain with at least three linearly independent positive solutions of equation  (\ref{Ground}) with $\Lambda = \Lambda_{\infty}$. Two of them are not  variational $\infty$--Ground states. 

  \paragraph{Convex Domains.} From now on we restrict ourselves to a \emph{convex} bounded domain. We need the crucial property that $v =\log(u)$ is concave. Therefore \emph{we assume that $u$ is a variational $\infty$--Ground state}.   By S. Sakaguchi's extension in \cite{S} of the Brascamp--Lieb Theorem the functions
  \begin{equation*}
    v_p\,=\,\log(u_p)\quad \text{and}\quad v\,=\,\log(u)
  \end{equation*}
  are concave,  where $u_p$ is a $p$-Ground state. (The concavity is preserved in the limit $p\to\infty$.)  The viscosity supersolutions $u > 0$ of the equation (\ref{Ground}) must satisfy the inequalities
  $$ \frac{|\nabla u|}{u}\,\geq\,\Lambda_{\infty}\quad\text{and}\quad \Delta_{\infty}u\,\leq\,0$$
  in the viscosity sense (i.\,e. for test functions touching $u$ from below). A  remarkable result of Y. Yu in \cite{Y}, Theorem 3.1, Lemma 3.5, is that if $u \in C^1(D)$ in an open set $D \subset \Omega$, then
  
  $$ \frac{|\nabla u(x)|}{u(x)}\,>\,\Lambda_{\infty}\quad \text{when}\quad x \in D$$
  and $\Delta_{\infty}u\,=\,0$ in $D$ (in the viscosity sense). In order to handle the singular set  $\{x\,|\,\,u\,\, \text{is not differentiable at}\,\, x\}$ one uses the operator
\[
  \Ss(x)\,=\, \underset{r \to 0}{\lim}\Bigl\{-\,\underset{y \in \partial B_r(x)}{\mathrm{min}}\frac{u(y)-u(x)}{r}\Bigr\}
  \]
  to define \emph{the contact set}
  $$\Upsilon\,=\,\bigl\{x\in \Omega\,|\,\,\Ss(x)\,=\,\Lambda_{\infty}u(x)\bigr\}.$$
According to Theorem 3.6 and Corollary 3.7 in [Yu], the contact set is closed and it has $N$-dimensional Lebesgue measure zero. Moreover, it contains the singular  set. At  points of differentiability we always have $S^-(x)\,=\,|\nabla u(x)|.$ In the open set $\Omega \setminus \Upsilon$, \, $\Delta_{\infty}u\,=\,0$ and $u \in C^1(\Omega \setminus \Upsilon).$

  The \emph{High Ridge} $H$ is defined as
  $$H\,=\,\bigl\{x \in \Omega\,|\,\mathrm{dist}(x,\partial \Omega)\,=\,R \bigr\}   ,$$
  where $R\,=\,\max \mathrm{dist}(x,\partial \Omega)$ is the radius of the largest ball that can be inscribed in $\Omega$. We shall use the
  $$\text{\sf{normalization}:}\quad  \underset{x \in \Omega}{\max}\,u\,=\,1.$$
  According to Theorem 2.4 in \cite{Y}, $H$ is also the set where $u$ attains its maximum:
  $$ H\,=\,\{x\in \Omega\,|\, u(x) \,=\,1\}.$$
The (ascending) streamlines of $u$   are defined as solutions of
  $$\frac{d\phantom{t}}{dt}\,\be(t)\,=\, \nabla u(\be(t))$$
    in $\Omega \setminus  \Upsilon$ but this does not work when a streamline enters $\Upsilon$, where the tangential direction $\nabla u$ is lost. Therefore we introduce \emph{fictitious streamlines} via a smoothing procedure  in Section \ref{sec:fic}. In this construction we use the function $v =\log(u)$ and not directly $u$.  A fictitious streamline goes from the boundary $\partial  \Omega$ towards the High Ridge $H$. In the plane it is unique, beginning as a proper streamline in  $\Omega \setminus   \Upsilon$ until it hits $\Upsilon$ at some point $x_{\Upsilon} = \al(t_{\Upsilon})$. After that it never leaves $\Upsilon$ on its way to $H$. They may meet and  join along a common arc  but cannot cross each others. There is always a boundary zone, say $ \mathrm{dist}(x,\partial \Omega)\,>\,\delta_0,$ where $\Delta_{\infty}u\,=\,0$ and $u \in C^1$, see page \pageref{zone}.

    \paragraph{Convex Polygons.} From now on we restrict our account to the plane. Suppose first that $\Omega$ is a convex polygon with corners $P_1,P_2,\dots,P_n$. The fictitious streamline
    $$\g_j\,=\,\g_j(t)\quad \text{from}\quad P_j\quad \text{to}\quad H, 
\qquad j = 1,2,...,n$$
    is called an \emph{attracting streamline}. The gradient $\nabla u$ exists on the boundary $\partial \Omega$, and  at the corners it vanishes. Having defined the  $\g_j$ with the aid of  $v$, we  use also   the streamlines of $u$ in the formulation of the  next theorem (outside $\Upsilon$ they are the same curves as those of $v$).

    \begin{thm}\label{polygon} The contact set has no points outside the attracting streamlines:
      $$\Upsilon\,\,\subset\,\, \g_1\cup \g_2\cup\cdots\cup \g_n$$
      Hence, $\Delta_{\infty}u\,=\,0$ and $u \in C^1$ outside the attracting streamlines. Moreover, a streamline $\be \,=\, \be(t)$ which is not one of the $\g_j$ cannot meet any other streamline before it meets a $\g_j$ (or reaches $H$). Its speed $|\nabla u(\be(t)|$ is constant until it joins $\g_j$.
    \end{thm}

    A similar theorem holds for a convex domain having the property that $|\nabla u|$ has only a finite number of maxima and minima along $\partial \Omega.$ See Theorem 3 in \cite{LL2} for the exact wording.
    
    The theorem has interesting consequences in convex polygons. Let $M_k$ denote the point on the side $P_kP_{k+1}$ at which $|\nabla u|$ attains 
its maximum:
    $$|\nabla u(M_k)|\,=\, \underset{x \in P_kP_{k+1}} {\max}|\nabla u(x)|.$$
    Now $|\nabla u|$ increases along the side  $P_kM_k$ and decreases along the side $M_kP_{k+1}$, see Lemma 18 in \cite{LL1}.  For lack of a better name, we call the streamlines starting at the boundary points $M_k$ for \emph{medians}. They have maximal speed.   At least two medians  are straight line segments joining the boundary to the High Ridge (because $u$ is squeezed between the distance function and a suitable cone function $|x-x_0|,\,\,x_0\in H$), but it is remarkable that they all start as straight lines:

    \begin{cor}[Medians] \label{median} The streamline starting at $M_k$ on $\partial \Omega$ is a straight line segment until it joins some   attracting streamline $\g_j$ or hits the High Ridge $H$.
    \end{cor}

    \begin{cor}[Arc length] \label{convex-concave} The arc of a streamline from the boundary $\partial \Omega$ till the first point $y$  at which it joins an attracting streamline or hits $H$, is either convex or concave. At the meeting point
      $$\frac{|\nabla u(y)|}{u(y)}\,=\,\frac{1}{S}$$
      where $S$ denotes the length of this arc. In particular, the length 
$S=1$ for all such arcs with the first meeting point $y \in \Upsilon$.
    \end{cor}

    A strange situation is possible. Suppose that an attracting streamline $\g_k$ hits $\Upsilon$ at the point $y_{\Upsilon}$. If  $y_{\Upsilon}$ does not belong to the Hige Ridge $H$, consider the arc of $\g_k$ between 
 $y_{\Upsilon}$ and $H$, which belongs to $\Upsilon$. All streamlines that \emph{first}
    hit $\g_k$ at this arc have the same length $S=1$ measured from $\partial \Omega$ to   $\g_k$.

    \paragraph{The Infinity Potential.} There is a related problem, which has been studied in \cite{L},\cite{LL1}, \cite{LL2}. The unique (see \cite{J}) solution $U$ of the boundary value problem
\[
      \begin{cases}
          \Delta_{\infty}U\,&=\,0  \quad\text{in}\quad \Omega \setminus 

H\\
       \phantom{\Delta_{a}}   U\,&=\,1 \quad\text{on}\quad H\\
       \phantom{\Delta_{a}}     U\,&=\,0 \quad\text{on}\quad \partial \Omega
        \end{cases}
   \]
      is called the $\infty$-potential. It has the advantage that $\nabla U$ is locally H\"{o}lder continuous in $\Omega \setminus H$ (cf. \cite{ES}) and its level curves are convex. (But it is not known whether $V = \log(U)$ is concave.) Using the  streamlines now defined through
      $$\frac{d}{dt}\,\be\,=\,\nabla U(\be),$$
      we have a similar result for $U$.

      \begin{thm} The counterpart to Theorem \ref{polygon} holds for the $\infty$-potential $U$ in a polygon.
      \end{thm}

      \medskip
      \emph{Proof:} This is Theorem 2 in \cite{LL2}.\qquad $\Box$

\begin{figure}
\begin{center}
\includegraphics[scale=0.11]{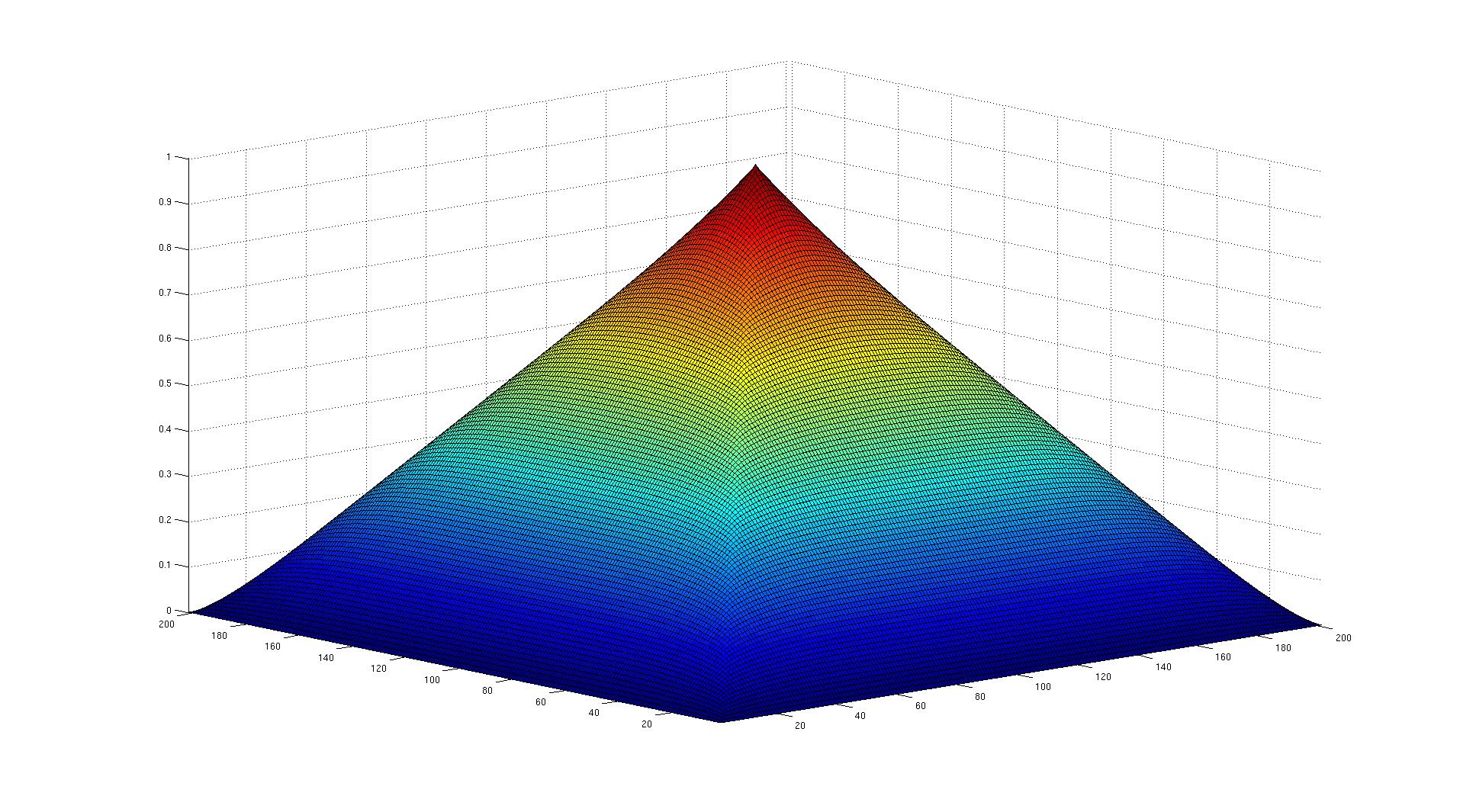}
\includegraphics[scale=0.032]{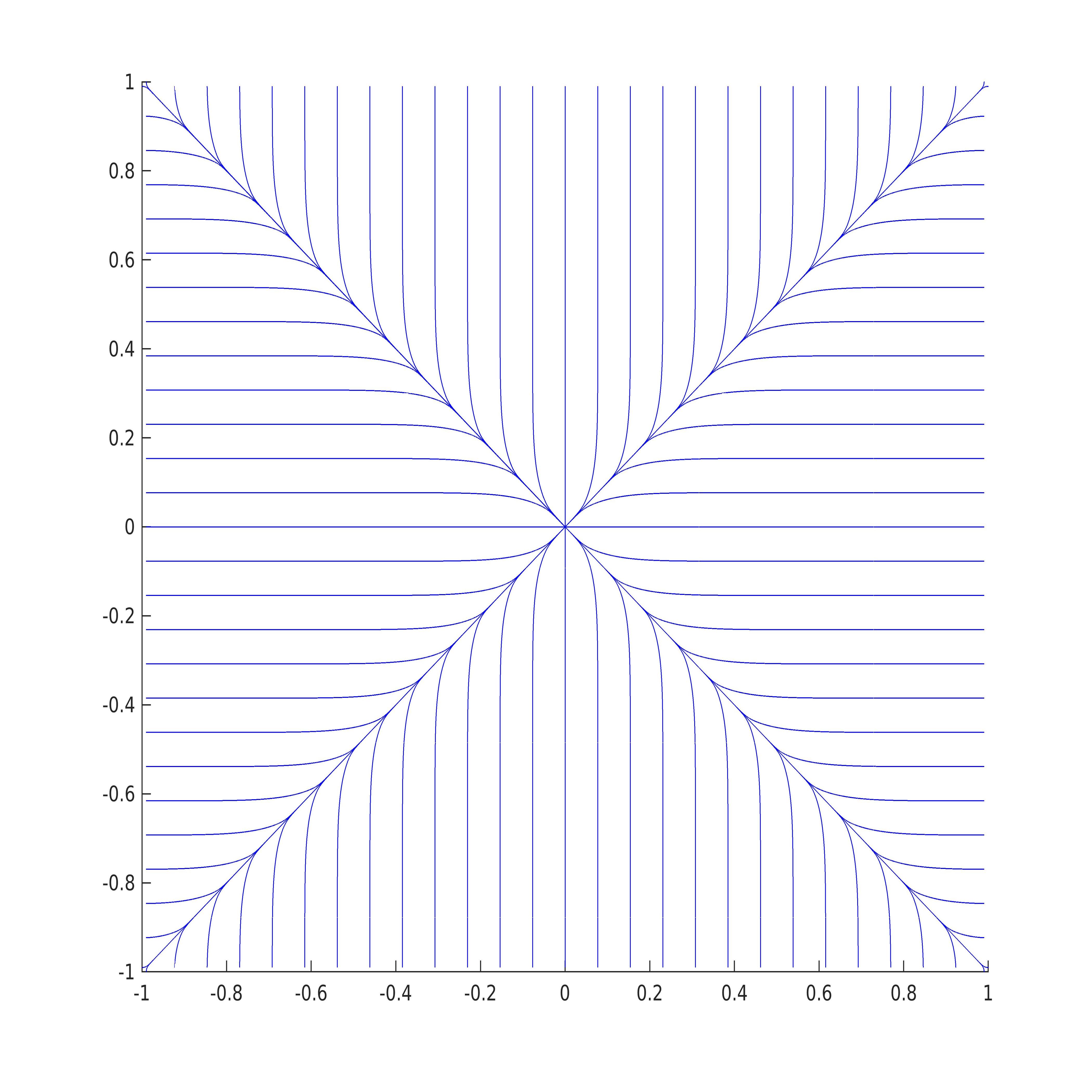}
\caption{The Infinity potential in the square and its streamlines.}\label{fig:square}
\end{center}
\end{figure}

      \bigskip
One may ask whether the functions coincide: Is $U=u$? This holds for stadium shaped domains. However, there are counterexamples even for convex polygons, see Theorem 3.3 in \cite{Y}. The problem is  intriguing  even for  a square, see Figure \ref{fig:square}. This case was our original impact. We cannot resist citing D. Hilbert: 
      \begin{quote}{\small Denn wer, ohne ein bestimmtes Problem vor Auge zu haben, nach Methoden sucht, dessen suchen ist meist vergeblich.}\end{quote}  

It should be kept in mind that the main difficulty is the lack of second partial derivatives of $u$ and $U$. We assume that the reader is familiar with the concept of viscosity solutions; see \cite{CIL} or \cite{K}.

\section{Preliminaries}

\paragraph{The $\Ss$--operator.} Let $\Omega$ denote a convex bounded domain in the Euclidean plane $\mathbb{R}^2.$ Let $u \in C(\overline{\Omega})\cap W^{1,\infty}_0(\Omega)$ be the normalized variational  $\infty$--Ground state. For the concave function $v=\log(u)$ we introduce the operator
$$\Ss(x)\,=\,\lim_{r\to 0}\Ss_r(x),\quad\text{where}\quad \Ss_r(x)\,=\, -\,\underset{x\in \partial B_{r}(x)}{\min}\, \frac{v(y)-v(x)}{r}.$$
In fact, $\Ss_r\,\searrow \, \Ss$ as $r\, \searrow 0.$ According\footnote{Y. Yu used the operator $S^-$ defined for $u$ in the place of our $\Ss$ defined for $v$. This has no bearing.} to [Yu] we have:
\begin{itemize}
\item $\Ss$ is continuous in $\Omega$ (Theorem 3.6).
\item At points of differentiability $\Ss(x)\,=\,|\nabla v(x)|.$
\item The contact set $\Upsilon\,=\,\bigl\{x\in\Omega\,|\,\Ss(x)\,=\,\Lambda_{\infty}\bigr\}$ is closed and has area zero (Corollary 3.7).
  \item The set where $u$ is not differentiable is contained in $\Upsilon$ (Lemma 3.5).
\end{itemize}
Also the High Ridge $H$ is contained in $\Upsilon$. It follows that
\begin{equation} \label{free} |\nabla v(x)|\,>\,\Lambda_{\infty}\quad \text{when}\quad x \in \Omega \setminus  \Upsilon.
\end{equation}

\paragraph{The  gradient estimate.} We shall employ the approximation
$$v_p = \log(u_p) \to  v = \log(u),\quad u_p \to u,\quad \|u_p\|_{\infty} \to \|u\|_{\infty} = 1,   $$
as $p\to\infty$, where $u_p\in C^{1,\kappa}(\overline{\Omega})$ is the $p$\,-\,Ground state in equation   (\ref{peqn}). By a general property for concave functions
$$\nabla v_p \rightarrow \nabla v\quad\text{locally uniformly in}\quad \Omega \setminus  \Upsilon.$$
In the next lemma, the bound is uniform in $p$, since
$$ \left(\lambda_p\,\mathrm{diam}(\Omega)\right)^{1/(p-1)} \|u_p\|_{\infty}\,\to\, \Lambda_{\infty}\|u\|_{\infty}.$$

\begin{lemma}\label{2021}. The estimate
  $$ |\nabla u_p(x)|\,\leq\,\left(\lambda_p\,\mathrm{diam}(\Omega)\right)^{\frac{1}{p-1}}\|u_p\|_{\infty}$$
  is valid when  $x \in \Omega$.
\end{lemma}

\emph{Proof:}  By general theory  (at least locally),  $\max|\nabla u_p|\,<\,\infty$ for $p$ fixed, but we need a bound that is stable as $p\,\to\,\infty$. To be on the safe side, we provide a simple proof using the convexity of the level curves. Given a point $z$ on the level curve $u_p(x) = c$, we shall construct a suitable function $w_p$ such that
$$ c = u_p(z)=w_p(z) \quad\text{and}\quad  c\leq u_p\leq w_p\quad \text{when} \quad u_p\geq c$$
from which we conclude by comparing difference quotients that $|\nabla u_p(z)|\leq |\nabla  w_p(z)|$. To construct a suitable $w_p$ we use the following version of  a  comparison principle.

If $u_p,\,w_p \in W^{1,p}(\Omega)$ are weak solutions in $D$ of
$$\Delta_pu_p\,+\,\lambda_p u_p^{p-1}\,=\,0,\quad \Delta_pw_p\,+\, \gamma\,=\,0,$$
and if $(u_p-w_p)^+\in W^{1,p}_0(\Omega)$ and $\lambda_p|u_p|^{p-1}\leq \gamma$, then we have $u_p\leq w_p$ in $D$. This see this, use  $(u_p-w_p)^+$ as test function in both equations and subtract these to obtain
\begin{gather*}
  \int_D\langle |\nabla u_p|^{p-2}\nabla u_p\,-\,|\nabla w_p|^{p-2}\nabla w_p,\nabla (u_p-w_p)^+ \rangle \,dx\\
  =\, \int_D(\lambda_p |u_p|^{p-2}u_p-\gamma)(u_p-w_p)^+\,dx\,\leq \, 0.
\end{gather*}
Then $w_p\geq u_p$ follows.

A solution $w_p$ that depends only on one coordinate, say $w_p(x) =  w_p(x_1)$, is easy to find. The equation reads $(p-1)|w_p'|^{p-2}w_p'' + \gamma = 0,$ where $'$ indicates $d/dx_1$. Upon solving the equation we obtain
$$w_p(x)\,=\,\frac{p-1}{p}\,\gamma^{\frac{1}{p-1}}\bigl\{\mathrm{diam}(\Omega)^{\frac{p}{p-1}}-(\mathrm{diam}(\Omega)\! -\! (x_1\!-\!z_1))^{\frac{p}{p-1}}\bigr\}\,+\,c
$$
when $0 \leq x_1-z_1 < \mathrm{diam}(\Omega)$. 
Finally, we rotate the coordinate system at $z$ so that the new $x_1$-axis points in the direction of the normal of the convex level curve at $z$.  By comparison, $w_p \geq u_p$.  At the point $z$ we now have
$$|\nabla u_p|\, =\, \frac{\partial u_p}{\partial x_1}\, \leq \,\frac{dw_p}{dx_1}\, =\, \gamma^{\frac{1}{p-1}}\mathrm{diam}(\Omega)^{\frac{1}{p-1}}$$
and to conclude, we choose $\gamma$ equal to $\lambda_p \|u_p\|_{\infty}^{p-1}.$ \qquad  $\Box$

\bigskip

To derive a lower bound we need  the concavity. We must avoid the points where $\nabla v_p  \approx 0$.

\begin{lemma}\label{1/c}
  Let $0<c<1$. There is a $p(c)$ such that the estimate
\[
    |\nabla v_p(x)|\,\geq\,\frac{1}{2\,\mathrm{diam}(\Omega)}\,\log\Bigl(\frac{1}{c}\Bigr),\qquad p >p(c),
  \]
  holds in the domain $v_p(x)\leq \log(c).$
\end{lemma}

\emph{Proof:} Fix an arbitrary point $x^*$ on the level curve
$v_p(x) = \log(c)$. 
Recall the normalization and let $y^*$ (independent of $p$) be some point of $H$ at which $v(y^*) = 0$. 
Define the unit vector
$$ e\,=\,\frac{y^*-x^*}{|y^*-x^*|}.$$
 The concavity of the  function $v_p(x^*+te)$ of $t$ implies
\begin{align*}
  \langle \nabla v_p(x^*),e\rangle\,=\, D_ev^{\e}(x^*)\,&\geq \frac{v_p(y^*)-\log(c)}{|y^*-x^*|}\\
  &\geq \frac{v_p(y^*)-\log(c)}{\mathrm{diam}(\Omega)}.
  \end{align*}
As $p \to \infty,\,\,v_p(y^*) \to  v(y^*) = \log(1) = 0$. Hence
$$v_p(y^*)\,>\,\frac{1}{2}\log(c)\quad\text{when}\quad p > p(c).$$
It follows that
$$|\nabla v_p(x^*)|\,\geq\,|\langle\nabla v_p(x^*),e\rangle|\,\geq\,\frac{1}{2\,\mathrm{diam}(\Omega)}\log(\frac{1}{c})$$
when $p > p(c)$ and $v_p(x^*) = \log(c)$. The lower bound is monotone in $c$ and so, choosing other auxiliary points $x^*$, we get an estimate valid when $v_p \leq \log(c)$. \qquad $\Box$

\begin{cor}\label{zone} Let $0<c<1$. The estimate
  $$\Ss(x)\,\geq\,\frac{1}{2\,\mathrm{diam}(\Omega)}\, \log\Bigl(\frac{1}{c}\Bigr)$$
  holds in the domain where $v(x)\,\leq\,\log(c)$.
\end{cor}

\emph{Proof:} In the set where $\nabla v(x)$ exists\, $|\nabla v_p(x)| \to |\nabla v(x)| \,=\Ss(x)$ even locally uniformly. Thus the the desired estimate  holds a.e., since the set of points of non-differentiability have  zero area. By the continuity of the $\Ss$--operator, the estimate is valid at \emph{every} point in question.
$\Box$

\medskip
This has an interesting consequence. There is a level $c_0$ such that
$$\Ss(x)\,>\, \Lambda_{\infty}\quad \text{when}\quad v(x)\,<\,\log(c_0).$$
Thus there is a \emph{boundary zone}, say the points $x \in \Omega$ with  $0<\mathrm{dist}(x,\partial\Omega)<\delta_0$ for some $\delta_0$, in which $u$ is a solution of $\Delta_{\infty}u\,=\,0$. In particular, $\nabla u$ is continuous up to the boundary $\partial \Omega$ there. See Theorem 1.1 in \cite{WY} and, for the corners,  Lemma 2 and Theorem 2 in \cite{HL} and also Theorem 7.1 in \cite{MPS}.
This boundary zone is  essential for  the procedure with quadrilaterals in Section \ref{Contact}.

 From Lemma \ref{1/c} it follows that $\nabla u_p\neq 0$ when $u_p < 1$. By standard regularity theory $u_p \in C^2$ when  $\nabla u_p\neq 0$. See Satz 10.3.2 in  \cite{Jo}   or, for an appealing proof in two variables, Theorem 12.5 on page 305 in \cite{GT}. The next proposition is later  needed to derive the Quadrilateral Rule.

\begin{prop} \label{Gauss} Suppose that the subdomain $D\subset \subset \Omega \setminus \Upsilon$ has a Lipschitz boundary $\partial D$ and let $\mathbf{n}$ denote the outer normal. Then
  $$\oint_{\partial D}|\nabla u|^{m-2}\langle \nabla u, \mathbf{n}\rangle\,ds \,\leq\,0$$
 is valid when  when $m\geq 2.$
\end{prop}
\bigskip
\emph{Proof:} This is essentially Proposition 1 in \cite{LL1},  but now the proof is much simpler, because
$$\nabla u_p\,\rightarrow\,\nabla u\quad\text{uniformly in}\quad \overline{D}.$$
This convergence can be deduced  from the fact that
$$\nabla v_p\,\rightarrow\,\nabla v\quad \text{locally uniformly in}\quad \Omega \setminus \Upsilon $$
according to a general theorem for concave functions.

By \cite{LMS}, \,  $\Delta_m u_p\, \leq \,0$ in the viscosity sense when $p\geq m\geq 2.$ Since $u_p \in C^2(D)$ for large $p$, this inequality holds pointwise in $\overline{D}$. It follows that we can use Gauss' Theorem and conclude that
$$\oint_{\partial D}|\nabla u_p|^{m-2}\langle \nabla u_p, \mathbf{n}\rangle\,ds \, = \, \int\!\!\!\!\int_{ D \,\,} \Delta_mu_p\, dxdy\,\leq\,0.$$
By the uniform convergence we can pass to the limit under the first integral sign to verify the desired inequality for $u$.\qquad $\Box$

\section{Fictitious Streamlines} \label{sec:fic}

A streamline $\al = \al(t)$ of $v$ should obey the rule
\begin{equation}\label{expected}\frac{d\al}{dt}\,=\,\nabla v(\al).\end{equation}
This is no problem in the boundary zone, but when $\al$ hits a point where $v$ is not differentable, the direction for the tangent is lost. It is convenient to prolong $\al$ along a \emph{fictitious streamline} the whole way to the High Ridge.
We shall construct it using the approximations $v_p = \log(u_p)$. Let $\al_p$ be a streamline for $v_p$. We have
\[
  \frac{d\al_p}{dt}\,=\, \nabla v_p(\al_p),\quad \al_p(t_0)\,=\,x_0.
\]
By the Picard--Lindel\"{o}f Theorem, the $\al_p$ is unique and approaches asymptotically  a point where $\nabla v_p = 0$. By direct calculation
$$\frac{d}{dt} v_p(\al_p(t))\,=\,\Big\langle \nabla v_p(\al_p(t)),\frac{d}{dt}\al_p(t)\Bigr\rangle\,=\,\bigl|\nabla v_p(\al_p(t))\bigr|^2,$$
and, since $v_p$ is a concave function,
$$\frac{d^2}{dt^2} v_p(\al_p(t))\,=\, \frac{d}{dt} \bigl|\nabla v_p(\,\al_p(t))\bigr|^2\,=\, 2\,\Delta_{\infty}v_p(\al_p(t))\,\leq0.$$ Thus $v_p(\al_p(t))$ is a concave function of $t$ and
\begin{equation}\label{concspeed}
  \bigl|\nabla v_p(\,\al_p(t_1))\bigr|^2\,\geq \bigl|\nabla v_p(\,\al_p(t_2)\bigr|^2, \quad t_2 > t_1.
\end{equation}

In order to take the limit as $p \to \infty$, we show that the family $\{v_p(\al_p)\}$ is locally equicontinuous. For
$t_0\leq t_1 <t_2< \infty$
\begin{align*}
  |\al_p(t_2) - \al_p(t_1)|\,=\,\Big\vert\int_{t_1}^{t_2}\nabla v_p(\al_p(t))\,dt\Big\vert
  \leq \int_{t_1}^{t_2}\big\vert\nabla v_p(\al_p(t))\big\vert\,dt \\ \leq\,(t_2-t_1)\max\big\vert\nabla v_p(\al_p(t))\big\vert
  \leq (t_2-t_1)   \big\vert\nabla v_p(\al_p(t_0))\big\vert   \leq C_{t_0}(t_2-t_1),                         
\end{align*}
where the uniform gradient bound $C_{t_0}$   comes from Lemma (\ref{2021}).
By Ascoli's Theorem we can extract a subsequence converging locally uniformly to some curve
$$\al(t)\,=\,\lim_{j \to \infty}\al_{p_{j}}(t), \quad \al(t_0) = x_0.$$
This is to be our fictitious streamline. Note that the limit $v(\al(t))\,=\,\lim v_p(\al_p(t))$
is a concave function in $t$.

Inside the open set $\Omega \setminus \Upsilon$ the constructed fictitious streamline obeys equation (\ref{expected}). To see this, notice that here $\nabla v_p$ converges locally uniformly to $\nabla v$. So does $\nabla v_p(\al_p(\tau)).$   We can take the limit in
$$\al_p(t)\,=\,\int^t\!\nabla v_p(\al_p(\tau))\,d\tau$$
and then differentate to obtain the desired equation. In particular, the tangent  $\tfrac{d}{dt}\al$ is continuous away from $\Upsilon$.

\paragraph{Stability.} The initial value problem is stable, if we use \emph{the same} $p_j$-subsequence for constructing two such curves, say
$$\al = \lim \al_p, \,\,\,\,\al(t_0) = x_0,\quad \muu = \lim \muu_p,\,\,\,\,\muu(t_0) = y_0.$$
The concavity of $v_p(x)$ implies that the expression
$$\frac{d}{dt}\,\bigl\vert\al_p(t)-\muu_p(t)\big\vert^2\,=\,2\,\big\langle\al_p(t)-\muu_p(t),
\nabla v_p(\al_p(t))-\nabla v_p(\muu_p(t))\big\rangle$$
is negative. It follows that
$$\bigl\vert\al(t)-\muu(t)\big\vert\,\leq\,|x_0-y_0|,\quad t> t_0.$$ 

Here we used the same subsequence, but 
it is well known that the solution of the initial value problem
$$
\frac{d\al}{dt}\,=\, \nabla v(\al),\quad \al(t_0)\,=\,x_0
$$
is unique, when $v$ is $C^1$ and concave. Provided that we stay away from $\Upsilon$ we can directly use the differential equation: If $\al$ and $\muu$ are two solutions with the same initial point $x_0$, the concavity again  implies that the expression
$$\frac{d}{dt}\,\bigl\vert\al(t)-\muu(t)\big\vert^2\,=\,2\,\big\langle\al(t)-\muu(t),
\nabla v(\al(t))-\nabla v(\muu(t))\big\rangle$$
is negative. It follows that
$$\al(t)=\muu(t),\quad t_{\Upsilon} > t> t_0.$$

In conclusion, the fictitious streamline $\al$ is unique  and independent of the used $\{p_j\}$--subsequence as long as it does not hit the contact set $\Upsilon$. Actually, it can never leave  $\Upsilon$ upon hitting it. 

We shall  later learn  that even after hitting $\Upsilon$, a fictitious streamline is unique, see page \pageref{unique}.

\begin{lemma}\label{caught} If\, $\al(t^*) \in \Upsilon$ for some $t^*$, then  $\al(t) \in \Upsilon$
  for all $t \geq t^*$.
\end{lemma}

\medskip
\emph{Proof:} If there is some point $\al(t_2)$ in $\Omega \setminus \Upsilon$  with $t_2 > t^*$, we follow $\al$ in its reverse direction till it hits $\Upsilon$. Put
$$t_1\,=\,\max\bigl\{t\,|\,t < t_2,\,\,\al(t)\in \Upsilon \bigr\}$$
In any case, $t^* \leq t_1 < t_2$. According to (\ref{concspeed})
$$ |\nabla v_p(\al_p(t_2))|\,\leq \, |\nabla v_p(\al_p(t))|\quad\text{when}\quad t_1<t<t_2.$$
Taking the limit as $ t \to t_1$   we get
$$ |\nabla v(\al(t_2))|\,\leq \, |\nabla v(\al(t))|\,=\,\Ss(\al(t))\, \to \,\Ss(\al(t_1))\,=\,\Lambda_{\infty}$$
by the continuity of $\Ss$.
By (\ref{free})  $|\nabla v(\al(t_2))|\,>\, \Lambda_{\infty}$. This is a contradiction.\qquad $\Box$

\bigskip
\paragraph{Intersecting the level curves.} Next we shall verify that a fictitious streamline, starting at a level low enough, reaches all higher level curves $v = \log(c)$ with $0<c<1$. Let $t_p(c)$ denote the unique time at which $v_p(\al_p(t))$ reaches the level  $\log(c)$. In other words  $v_p\bigl(\al_p(t_p(c))\bigr) \,=\, \log(c).$ For two levels $0<c_1<c_2<1$ we have
$$\int \limits _{t_p(c_1)}^{t_p(c_2)}|\nabla v_p(\al_p(t))|\,dt\,\geq\,\bigl(t_p(c_2)-
t_p(c_1)\bigr)\frac{\log(\frac{1}{c_2})}{2\,\mathrm{diam}(\Omega)}.$$
 by Lemma \ref{1/c}.  On the other hand
\begin{align*}
  \int \limits_{t_p(c_1)}^{t_p(c_2)}|\nabla v_p(\al_p(t))|\,dt\,&\leq \bigl(t_p(c_2)-
  t_p(c_1)\bigr)^{\frac{1}{2}}\Bigl(\int \limits_{t_p(c_1)}^{t_p(c_2)}|\nabla v_p(\al_p(t))|^2\,dt\Bigr) ^{\frac{1}{2}}\\
 &= \bigl(t_p(c_2)-\,  t_p(c_1)\bigr)^{\frac{1}{2}}\Bigl(v_p(\al(t_p(c_2)))-
 v_p(\al(t_p(c_1)))\Bigr)^{\frac{1}{2}}\\
 &=\Bigl( \bigl(t_p(c_2)- t_p(c_1)\bigr)\log(\frac{c_2}{c_1})\Bigr)^{\frac{1}{2}}.
\end{align*}
Combining the above estimates we see that
$$\bigl\vert t_p(c_2)- t_p(c_1)\bigr\vert\,\leq\biggl(\frac{2\,\mathrm{diam}(\Omega)}{\log(\frac{1}{c_2})}\bigg)^2\log(\frac{c_2}{c_1}).$$
We conclude that the family $\{t_p(c)\}$ is locally equicontinuous. It is also seen to be equibounded, if we let all curves start at the same point. Hence we can again use Ascoli's Theorem to extract a subsequence $p_1,p_2,p_3,\dots$  for which
$$ t(c)\,=\,\lim_{j\to \infty} t_{p_j}(c)$$
locally uniformly. This means that $v\bigl(\al(t(c))\bigr)\,=\,\log(c),$ whenever $0 < c < 1$. Thus we have reached all levels.

\emph{A fictitious streamline intersects every level curve at exactly one point}, because the function
$v(\al(t))$ is \emph{strictly} increasing in $t$. The strict inequality $ v(\al(t_2)) > v(\al(t_1))$, where $t_2 > t_1$,   follows by sending $p$ to $\infty$ in 
\begin{align*}
  v_p(\al_p(t_2))- v_p(\al_p(t_1))\,&=\,\int\limits_{t_1}^{t_2}|\nabla v_p(\al_p(t))|^2\,dt \\ &\geq\,(t_2-t_1)\Bigl(\frac{1}{2\,\mathrm{diam}(\Omega)}\log(\frac{1}{c})\Bigr)^{\frac{1}{2}},
  \end{align*}
valid when  $v(\al_p(t_2)) \leq \log(c).$

\paragraph{Passage through all points.} Through every point $x_0$ outside $H$ there passes a fictitious streamline going from $\partial \Omega$ to $H$. To see this, consider the streamlines $\al_p =\al_p(t)$ first starting at the point:
$\al_p(t_0) = x_0$, say outside the High Ridge $H$. They go to $H$. They can be prolonged in the reverse direction by solving
the equation 
$$\frac{d}{dt}v_p(\al_p)\,=\,- \nabla v_p(\al_p)$$
with a \emph{ minus} sign. Now $t_0$ becomes an interior point of the parameter interval. This is classical theory. As we know, $\al_p \to \al$ locally uniformly. In particular $\al(t_0) = x_0$ as desired. However, the arc \emph{below} $x_0$ is not unique!

\paragraph{Uniqueness.} \label{unique} Finally, let us prove that a fictitious streamline $\al$ with a given initial point is unique. Suppose that there are two such curves $\al_1$ and $\al_2$ with the same initial point.  We can assume that the initial point is in $\Omega \setminus \Upsilon$. Then we already know that they coincide at least up to a point $x_{\Upsilon} = \al_1(t_{\Upsilon})  = \al_2(t_{\Upsilon}) $ which belongs to $\Upsilon$. After this, the curves stay in $\Upsilon$ (Lemma \ref{caught}).     If the uniqueness is violated, there must be two different  points on the same level curve $\om$ so that $a = \al_1(t_a)$ and $b = \al_2(t_b)$. Follow both curves in the reverse direction from $\om$ till the  point $c$ where they  departed. In any case, $c \in \Upsilon$. Consider the curved triangle $a,b,c$ bounded by the arcs of $\om, \al_1,$ and $\al_2.$ Now the arcs of the two fictitious streamlines belong to $\Upsilon$.

We claim that  all points in this triangle belong to $\Upsilon$. Since the triangle has positive area, this contradicts the fact that $\mathrm{area}(\Upsilon) = 0.$  To this end, suppose that there is a point $y$ in the triangle  not belonging to $\Upsilon$. A fictitious streamline that passes through  $y$ and starts for instance in the boundary zone must cross $\al_1$ or $\al_2$ before reaching  $y$. That means that after the crossing this new fictitious streamline stays in  $\Upsilon$. In particular
we must have $y \in \Upsilon$. This shows that there was no bifurcation.

\paragraph{Summary.} We remark that  fictitious streamlines may meet and then continue along  a common arc. But they cannot cross. The following properies have been proved:
\begin{itemize}
\item  At each point in $\Omega\setminus H $ there starts a unique fictitious streamline
    $\al$.
\item $v(\al(t))$ is strictly increasing as a function of $t$.
 \item Through every point
    there passes a  fictitious streamline, not necessarily unique.
\item $\al$ intersects a level curve  only  once. (It can be prolonged so 

that it intersects every level curve exactly once.)
\item If $\al$ starts in $\Omega \setminus \Upsilon$, it first  hits $\Upsilon$ at a unique point  $\al(t_{\al})$ and the speed $|\nabla v(\al(t))|$ is decreasing as long as $t\leq t_{\al}$.
\item $\Ss(\al(t))\, =\,  \Lambda_{\infty}\quad\text{when}\quad t \geq t_{\al}.$
\end{itemize}
At the end, it will appear that the only fictitious streamlines are the attracting ones: $\g_1,...,\g_n.$ (It may happen that $\nabla v$ does not exist on the curve.)

\section{ Non-decreasing Speed $|\nabla u(\be(t))|$}

It is convenient to consider the streamlines of $u_p$ and $u$, especially near $\partial\Omega$ where   $v \approx - \infty$:
  $$\frac{d}{dt}\,\be_p(t)\,=\,\nabla u_p(\be_p(t)),\qquad \frac{d}{dt}\,\be(t)\,=\,\nabla u(\be(t)).$$
  Needless to say, $u$ and $v =\log(u)$ have the same streamlines in  $\Omega\setminus\Upsilon$, though with different parametrizations. Indeed, in this case
  $$\frac{d}{d\tau}\,\al(\tau)\,=\,\nabla v(\al(\tau)); \quad \be(t)\,=\,\al\bigl(\tau(t)\bigr),\quad \tau(t)\,=\,\int_{t_0}^t u(\al(\sigma))\,d\sigma.$$
  Recall that $u_p \in C^2$ where $\nabla u_p \neq 0.$

  \begin{thm}\label{increasing} Along  an arc of the streamline $\be = \be(t)$,  comprised in $\Omega \setminus \Upsilon$, the function $u(\be(t))$ of $t$ is convex and the speed is non-decreasing:
    $$|\nabla u(\be(t_2))|\,\geq\,|\nabla u(\be(t_1))|,\quad\text{when}\quad t_2>t_1.$$
  \end{thm}

  \medskip
  \emph{Proof:} Recall that $u_p \rightarrow u$ uniformly in $\Omega$ and that $\nabla u_p \rightarrow \nabla u$ locally uniformly in $\Omega \setminus \Upsilon$. Let $K \subset \subset \Omega \setminus \Upsilon$.   We have the pointwise bounds
  $$ \Lambda_{\infty}\,< \,\frac{|\nabla u(x)|}{u(x)}\, <\, C_K, \qquad x 
\in K.$$
 The uniform convergence in $K$ implies that for some suitable constant $\gamma > 0$
  $$ \Lambda_{\infty}  + \gamma\,< \, \frac{|\nabla u_p|}{u_p}\,<\,C_K\quad \text{in}\quad K\quad \text{when}\quad p > p(\gamma,K).$$
  It is enough to prove the theorem for a closed subarc of $\be$. We can select $K$ above so that taking $p$ large enough $\be_p(t) \in K$ when $t_1\leq t\leq t_2,$  say. Indeed, from
  $$\be_p(t')-\be_p(t'')\,=\,\int_{t'}^{t''}\! \nabla u_p(\be_p(t))\,dt$$
  it again follows with Ascoli's Theorem
  that $\be =\lim\limits_{p \to \infty}\be_p$ uniformly in $[t_1,t_2].$ Then 
  $$\be(t')-\be(t'')\,=\,\int_{t'}^{t''}\! \nabla u(\be(t))\,dt$$
  and
  $$\frac{d}{dt}\,\be\,=\,\nabla u(\be),$$
  as it should.
  
  We shall show that for some suitable  constants $\kappa_p\rightarrow 0$ as $p\rightarrow 0$ the function
  $$u_p\bigl(\be_p(t)\bigr) + \frac{1}{2} \kappa_pt^2,\quad  t_1\leq t\leq t_2, $$
  is convex in the variable  $t$. We can differentiate twice to see that
  \begin{align*}
    \frac{d}{dt}\,u_p(\be_p(t))\,&=\,\Bigl\langle \nabla u_p(\be_p(t)), \frac{d\be_p(t)}{dt}\Bigr\rangle\,=\,|\nabla u_p(\be_p(t))|^2,\\
    \frac{d^2}{dt^2}\,u_p(\be_p(t))\,&=\, \frac{d}{dt}\,|\nabla u_p(\be_p(t))|^2\,=\,   2\,\Delta_{\infty}u_p(\be_p(t)).
  \end{align*}
  Using the equation $\Delta_p u_p + \lambda_pu_p^{p-1} = 0$ we can write
  \begin{align*}
    \Delta_{\infty}u_p\,=&\,-\,\frac{|\nabla u_p|^2\Delta u_p}{p-2}-\frac{\lambda_pu_p^{p-1}}{(p-2)|\nabla u_p|^{p-4}}\\
    \geq& \,-\,\frac{\lambda_p}{p-2}|\nabla u_p|^3\frac{1}{|\nabla v_p|^{p-1}},
  \end{align*}
  because $\Delta u_p\,\leq\,0$ for $p\geq 2$ according to Corollary 3.14 in [LMS]. Since $\sqrt[p]{\lambda_p} \rightarrow  \Lambda_{\infty}$ we have
  $$\lambda_p\,\leq\,\Bigl( \Lambda_{\infty}+\frac{\gamma}{2}\Bigr)^p $$
  for large $p$. Finally,
  $$2\,\Delta_{\infty}u_p\,\geq\,-\frac{2C_K^3}{p-2}\biggl(\frac{\Lambda_{\infty}+\frac{\gamma}{2}}{\Lambda_{\infty}+\gamma}\biggr)^{p-1}\!\!\Bigl( \Lambda_{\infty}+\frac{\gamma}{2}\Bigr)\,\equiv\,-\kappa_p.$$
  It is clear that $\lim \kappa_p = 0$ as $p \to \infty$.
  
  Integrating, we see that
  $$\frac{d}{dt}\,\bigl(|\nabla u_p(\be_p(t))|^2 + \kappa_pt\bigr)\,\geq\,0$$
  and hence $u_p(\be_p(t)) +\tfrac{1}{2}\kappa_pt^2$ is convex. So is its limit $u(\be(t))$. The theorem follows.\qquad $\Box$
  
  \section{Reducing the Contact Set $\Upsilon$}\label{Contact}

 \paragraph{The Quadrilateral Rule.}  An important device in our method is the ``Quadrilateral Rule'', which we derived in [LL2] for solutions of the equation $\Delta_{\infty}V = 0$ in convex ring domains. Consider a quadrilateral $Q$ bounded by arcs of  two streamlines $\be$ and $\et$ and two level curves $\om$ and $\si$. 

\begin{figure}[h!]
  \begin{center}
 
	\begin{tikzpicture}[domain=-5:5,scale=1.2]
	       
\fill (3,0) circle[radius=2pt] node[below] {$e$};
\fill (-3,0) circle[radius=2pt] node[below] {$b$};
\fill (2,4) circle[radius=2pt] node[above] {$e'$};
\fill (-2,4) circle[radius=2pt] node[above] {$b'$};
\draw[color=red, thick] (-1.6,-0.46) to[out=85,in=-92] (-2.04,2.5);
\draw[color=red, thick] (0.5,-0.71) .. controls (0.8,1.5) .. (0.5,3.25); 
\draw (-3,0) .. controls (0,-1) .. (3,0);\fill (-1.2,-0.8) node[below] {$\si$};
\draw (-3,0) .. controls (-2,2) .. (-2,4);\fill (-2.2,2) node[left] {$\be$};
\draw (3,0) .. controls (3,2) .. (2,4);\fill (3,2) node[right] {$\et$};
\draw (2,4) .. controls (0.5,3) .. (-2,4); \fill (-0.5,3.6) node[above] {$\om$};





	\end{tikzpicture}
	
\end{center}
\caption{The quadrilateral $bee'b'$ with two possible streamlines (in red).}
	\end{figure}
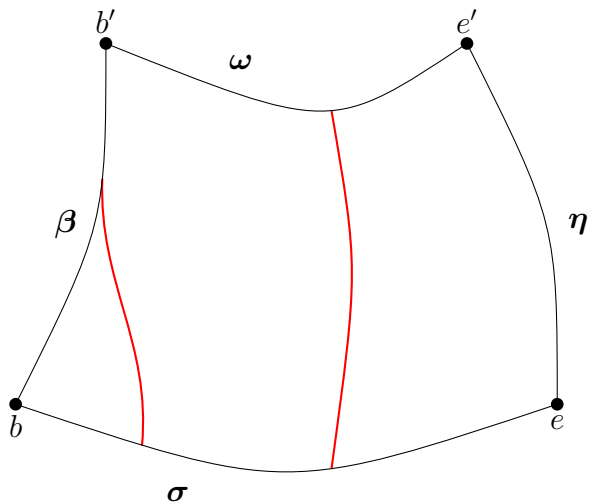
\phantom{a}
The proof required
\begin{itemize}
\item The Fundamental Lemma (Proposition \ref{Gauss}).
\item Increasing speed along streamlines (Theorem \ref{increasing}).
\item Monotonicity of $|\nabla u(\si)|$ along the lower arc ([LL2, Lemma 18]).
\end{itemize}
Thus the same proof as in [LL2] works well now for our function $u$, provided that $\nabla u$ is continuous also on the boundary curves of the quadrilateral, which means that $\Upsilon$ has to be avoided.  However, this rule is valid also if  the streamlines have common parts with  $\Upsilon$. Now Proposition \ref{Gauss} is not directly valid, but the only conclusion needed is that in a quadrilateral
$$\underset{b'e'}{\max}|\nabla u(\om)|\,\leq \underset{be}{ \max}|\nabla u(\si)|.$$
This can be proven by $p$-approximation, since integrals along $p$-streamlines disappear. We shall avoid this little addition here; this causes some extra concern in the choice of auxiliary quadrilaterals later.

We use the streamlines of $u$ below.

  \begin{lemma}[Quadrilateral Rule]\label{rule} Suppose that the streamlines $\be$ and $\et$ together with the level curves $\si$ (lower level) and $\om$ (upper level) form a quadrilateral $Q \subset \overline{\Omega}\setminus \Upsilon$ with vertices $b,e,e',b'$. Assume that $|\nabla u(\si)|$ is monotone on the arc $\overline{be}$ of $\si$. Then no streamline with initial point on the arc $\overline{be}$ (excluding  $b$ and $e$)  can meet any other streamline strictly inside the quadrilateral. Such a
    streamline  has constant speed till it joins  $\be$ or $\et$, or reaches $\om$.
  \end{lemma}

  \medskip
  We shall use an immediate consequence:
  
  \begin{cor} 
   The speed is monotonic along the arc $\overline{b'e'}$ on the  upper level curve $\om$.
    \end{cor}

  \medskip

  It is possible that the quadrilateral degenerates to a triangle bounded by $\si, \et$ and $\be$, when the upper level curve is reduced to a point: $b'=e'$. In this situation we have the corresponding \textbf{Triangular Rule}. Its proof follows from an approximation by quadrilaterals with upper level arcs slightly below the level of the corner $b'=e'.$
   
  \paragraph{The Polygon.}  Let us now apply the Quadrilateral Rule to the case when the domain $\Omega$ is a convex polygon with corners $P_1,P_2,\dots P_n$ and set $P_{n+1} = P_1$. Recall that we have a boundary zone, where
   $u$ is $C^1$ up to the boundary $\partial \Omega$ and $\Delta_{\infty}u\,=\,0$ there. In particular, the gradient $\nabla u$ is continuous along any side, say $P_1P_2$ and at the corners  $|\nabla u(P_j)| = 0$. 
   Now $|\nabla u|$ has a maximum along the side   $P_1P_2$, say at the point $M$.   By\footnote{The proof only required the convex level curves.} [LL2, Lemma 18], the normal derivative
   $$\frac{\partial u}{\partial n}\,=\, |\nabla u|$$
   is monotone along the segments $P_1M$ and $MP_2$. The streamline $\g$ starting at $P_1$ is called an attracting streamline. With the help of $v$ it can  be prolonged as a fictitious streamline
   going  the whole way to $H$. In the same way, a (fictitious) streamline $\muu$  starts at $M$ and approaches $H$ asymptotically. We consider the ``sector'' $D$  bounded by the two streamlines, the edge $P_1M$ and possibly $H$, if they do not join.

   \begin{lemma}\label{sweep} If there is a point $z\in D$ that belongs to $\Upsilon$ but does not lie on the attracting streamline $\g$, then the whole level arc $\{x\in D \,|\,u(x) = u(z)\}$  from $z$ to $\g$ belongs to $ \Upsilon$.
   \end{lemma}

   \medskip
   \emph{Proof:} Suppose, on the contrary, that there were on the level arc a point $y$   at which $|\nabla(v(y))| > \Lambda_{\infty}$. The streamline $\be$ from the side $\overline{P_1M}$ to $y$ cannot contain any points of $\Upsilon$, because that would imply $y\in \Upsilon$. Recall that $ \Upsilon$ has positive distance to the boundary. Let $\om$ be the highest level arc connecting   $\muu$ and $\be$ such that there are no points of $\Upsilon$ in the \emph{open} quadrilateral bounded by $\muu, \be, \om$ and $P_1M$. There is a point $c$ on $\om$ that lies in $\Upsilon$, since $\Upsilon$ is closed.  (It is possible, for example,  that $c = z$.). Let $b$ be the point where $\be$ and $\om$  intersect. 

\begin{figure}[h!]
  \begin{center}
 
	\begin{tikzpicture}[domain=-5:5,scale=2.]
	       
\fill (3,0.5) circle[radius=2pt]; \fill (3,0.4) node[below] {$M$};
\fill (-1.5,0.12) circle[radius=2pt]; \fill (-1.5,0)  node[below] {$P_1$};
\draw[color=red, thick] (0.35,0.28) to[out=90,in=-85] (0.3,2.4);

 \fill (0.3,1.5) node[left] {$\be$};

\draw[color=red, thick] (0.3,2.4) to[out=95,in=-88] (0.25,3.22); 
\fill (0.25,3.21) circle[radius=2pt];\fill (0.35,3.22) node[above] {$y$};

\fill (1.2,3.48) circle[radius=2pt];\fill (1.3,3.55) node[above] {$z$};

\draw[color=red, thick] (1.5,0.37) to[out=91,in=-72] (1.2,3.04); 

\draw (-1.5,0.12) -- (3,0.5);
\draw (-1.5,0.12) .. controls (-1.1,2) .. (-1,3.6);\fill (-1.35,2) node[left] {$\g$};
\draw (3,0.5) .. controls (3,2) .. (2,4);\fill (3,2) node[right] {$\muu$};
\draw (2,4) .. controls (0.5,3) .. (-1,3.6); 
\draw[color=blue] (0.29,2.38) .. controls (1.15,2.4) .. (2.44,3.1); \fill (1.95, 2.75) node[below] {{\color{blue}$\om_\e$}};
\draw[color=blue] (0.27,2.85) .. controls (1,2.9) .. (2.21,3.6); \fill (1.8, 3.05) node[above] {{\color{blue}$\om$}};


\fill (0.27,2.85) circle[radius=2pt]; \fill (0.22,2.85)  node[left] {$b$};
\fill ((0.3,2.4) circle[radius=2pt]; \fill (0.26,2.4) node[left] {$b_\e$};
\fill (1.34,2.54) circle[radius=2pt];
\fill (1.41,2.6) node[above] {$c_\e$};

\fill (1.2,3.05) circle[radius=2pt];
\fill (1.25,3.1) node[above] {$c$};





	\end{tikzpicture}
	
\end{center} 
\caption{The curves giving rise to the quadrilateral bounded by $\muu, \be$, $P_1M$ and $\om_\e$.}
  \label{fig:quadc}  
	\end{figure}
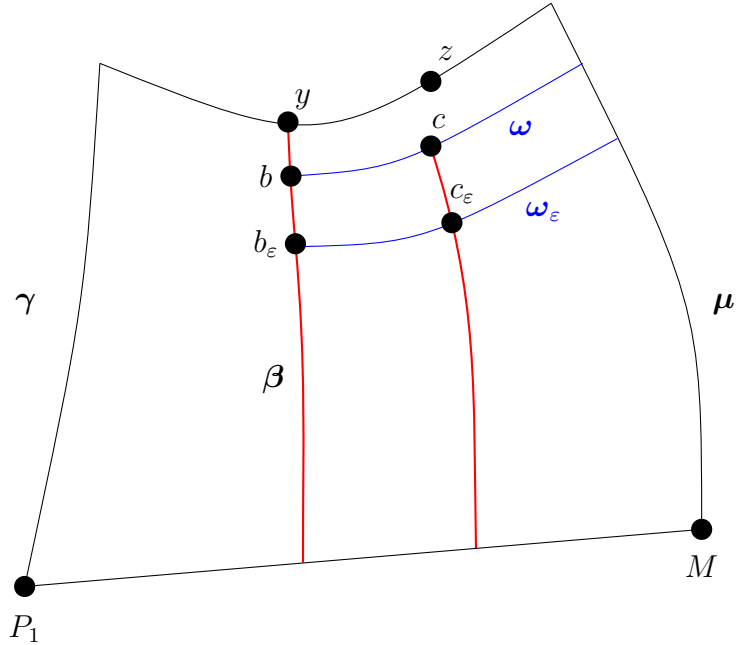

   The advantage of the arrangement is that we can now use the Quadrilateral Rule in the quadrilateral bounded again by  $\muu, \be$, and   $P_1M$ but with $\om$ now replaced by a slightly lower level curve  $\om_{\e}$, say  $\{u(x) =u(c)-\e\}$. Let $c_{\e}$ be the point where the streamline from $P_1M$ to $c$ intersects $\om_{\e}$. At the point  $b_{\e}$ where this new level curve intersects $\be$ we
   must have
   $$|\nabla u(b_{\e})|\,\leq\,|\nabla u(c_{\e})|,\qquad u(b_{\e}) = u(c_{\e}), $$
   by the remark after the Quadrilateral Rule (Lemma \ref{rule}).
   Thus
   $\Ss(b_{\e})\,\leq\,\Ss(c_{\e})$. Letting $\e \to 0$ we get by the continuity of the $\Ss$--operator
   $$\Ss(b)\,\leq\,\Ss(c)\,=\,\Lambda_{\infty}$$
   Thus $\Ss(b) = \Lambda_{\infty}$. Thus we arrive at the contradiction $b\in \Upsilon$. (The point $b$ was on $\be$, at a level that had not reached $\Upsilon$.) Thus the point $y$ does not exist. This proves that the whole arc is in the contact set $\Upsilon$. \qquad $\Box$

   \bigskip

   \begin{prop}[Structure of $\Upsilon$]\label{last}
   The contact set $\Upsilon$ is a subset of the union  $\g_1\cup\g_2\cup \cdots\cup \g_n $ of the attracting streamlines.
   \end{prop}

   \medskip
   \emph{Proof:} Consider again the ``sector'' $D$ bounded by the streamlines $\g_1,\,\muu$ and the
   side $P_1M$. Suppose that there exists a point $z \in \Upsilon$ that does not lie on $\g_1$. We shall show that this is impossible.  Let $\si$ denote the level arc from $z$ to $\g_1$.   Then Lemma \ref{sweep} implies that  $\si$ belongs to $ \Upsilon$. We shall see that this forces the contact set to have positive area.

   Take two points on $\si$, say $a_1$ and $a_2$. Construct the fictitious streamlines $\al_1$ and and $\al_2$ joining these points with $H$. They cannot meet on the level curve $\si$, since they cross a level only once. Follow the streamlines til they join each other or reach $H$. All points in the ``triangle'' bounded by $\si,\,\al_1,\,\al_2$ and possibly $H$ must belong to $\Upsilon$.
   Indeed, otherwise there would be an open component in the triangle where $|\nabla v| > \Lambda_{\infty}$. A streamline in that component travelled in the reverse direction cannot escape: it must hit a point in $\Upsilon$. But this contradicts Lemma \ref{caught}, according to which a streamline cannot leave  $\Upsilon$.  

   We have reached a contradiction, since the area of the triangle is $> 0$ but it is contained in $\Upsilon$. Recall that $\mathrm{area}(\Upsilon) = 0.$\qquad $\Box$

   \paragraph{Proof of Theorem 1.}  The first part of the theorem is included in Proposition \ref{last}. Also, we know that $u$ is of class $C_{loc}^{1,\kappa}$ outside $\Upsilon$ and that the equation $\Delta_{\infty} u\,=\,0$ is valid there.
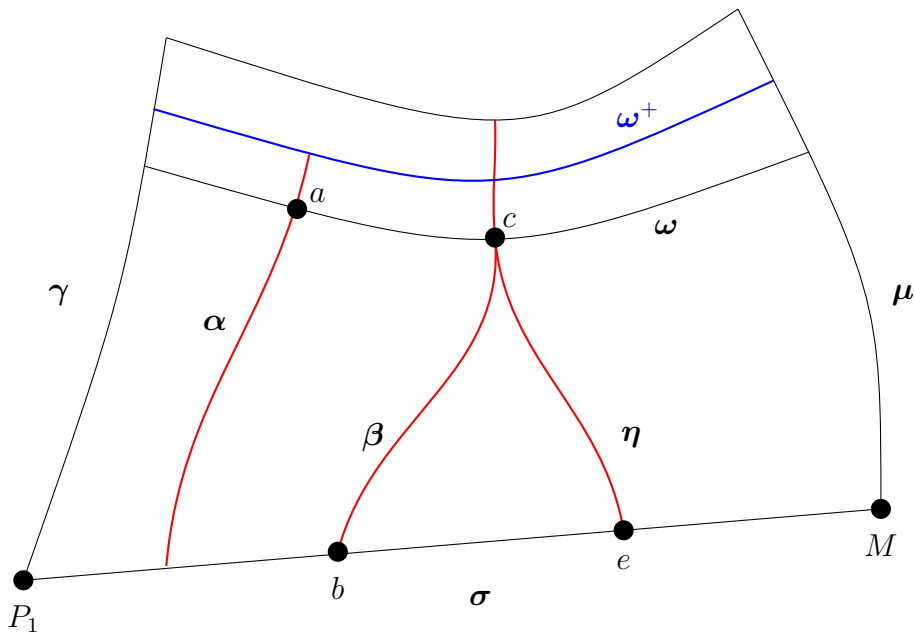
\begin{figure}[h!]
  \begin{center}
 
	\begin{tikzpicture}[domain=-5:5,scale=1.9]
	       
\fill (3,0.5) circle[radius=2pt]; \fill (3,0.4) node[below] {$M$};
\fill (-3,0) circle[radius=2pt]; \fill (-3,-0.1)  node[below] {$P_1$};
\draw[color=red, thick] (-0.8,0.2) to[out=75,in=-85] (0.3,2.4); \fill (-0.4,1) node[left] {$\be$};
\draw[color=red, thick] (0.3,2.4) to[out=95,in=-88] (0.3,3.22); 

\draw[color=red, thick] (-2,0.1) to[out=85,in=-102] (-1,2.98); \fill (-1.5,1.8) node[left] {$\al$};
\draw[color=red, thick] (1.2,0.35)  to[out=100,in=-85] (0.3,2.4); \fill (1.4,1) node[left] {$\et$};
\draw (-3,0) -- (3,0.5);\fill (0.2,0) node[below] {$\si$};
\draw (-3,0) .. controls (-2.3,2) .. (-2,3.8);\fill (-2.6,2) node[left] {$\g$};
\draw (3,0.5) .. controls (3,2) .. (2,4);\fill (3,2) node[right] {$\muu$};
\draw (2,4) .. controls (0.5,3) .. (-2,3.8); 
\draw (-2.15,2.9) .. controls (0.3,2.2) .. (2.5,3); \fill (1.5,2.6) node[below] {$\om$};
\draw[color=blue, thick] (-2.09,3.3) .. controls (0.3,2.6) .. (2.25,3.5); \fill (1.3,3.1) node[above] {{\color{blue}$\om^+$}};

\fill (-0.8,0.2) circle[radius=2pt]; \fill (-0.8,0.1)  node[below] {$b$};
\fill (1.2,0.35) circle[radius=2pt]; \fill (1.2,0.25) node[below] {$e$};
\fill ((0.3,2.4) circle[radius=2pt]; \fill (0.4,2.4) node[above] {$c$};
\fill (-1.085,2.6) circle[radius=2pt];
\fill (-1.075,2.7) node[right] {$a$};





	\end{tikzpicture}
	
\end{center} 
\caption{The curves giving rise to the quadrilateral with $c$ as an interior point.}
  \label{fig:quadc}  
	\end{figure}
   It remains to address the streamlines of $u$. First we show that there are no meeting points outside $\Upsilon$. On the contrary, suppose that two streamlines $\be$ and $\et$ with initial points $b$ and $e$ on the side $P_1M$ of the polygon \emph{first} meet at a point $c$ inside $\Omega  \setminus  \Upsilon$. We shall construct a quadrilateral in which $c$ is strictly inside. Then the Quadrilateral Rule  shows that $c$ does not exist.  To this end,   let $\g$ be the streamline starting at the corner $P_1$. Assume first that the meeting point $c$ is not on  the streamline $\muu$ with  initial point $M$.  Let $\om$ denote the level curve through $c$. Select a point $a$ on $\om$ strictly between $\g$ and $c$ and let $\al$ denote the streamline through $a$. By the uniqueness of streamlines outside $\Upsilon$, $\al$ cannot meet $\be$ below $\om$. The point $a$ has a neighborhood that is disjoint from both $\g$ and $\be$ and $\et$. Thus we can find a level curve $\om^+$ above $\om$ such that $\al$ does neither meet $\be$ nor $\g$ before it has reached $\om^+$. 
   Now the desired quadrilateral with $c$ as an interior point is the one bounded by  $P_1M$\, ($\si$ on the picture), $\muu$, $\om^+$, and $\al$. See Figure \ref{fig:quadc}.

   Moreover, we can now conclude from the Quadrilateral Rule, that the speed $|\nabla u(\be)|$ of any streamline is constant until the streamline joins $\g$ or $\muu$.  It remains to exclude $\muu$. So, if the streamline $\et$  starting at $e$ on the side $P_1M$ first meets $\muu$ at a point $m$, not on $\Upsilon$, both must have the same speed: $|\nabla u(\et)|\,=\,|\nabla u(\muu)|$ at the point $m$. It follows that both have the same  constant speed up to the meeting point. By the monotonicity, we can conclude that $|\nabla u|$ is constant along the segment between $e$ and $M$ on the side $P_1M$. Now it follows that each streamline starting on this segment has the same constant speed in the triangle bounded by $eM$, $\muu$, and $\et$. Therefore $u$ is a solution to the Eikonal Equation $|\nabla u|\,=\,C$ (a constant) in the triangle. Since $u$ is of class $C^1$ here, the streamlines must be segments of straight lines inside the triangle. This is a well-known property of the Eikonal Equation, see Lemma 1 in \cite{A2}. But this is impossible if $u$ is differentiable at the meeting point, since the tangents  have different directions. In conclusion, the only possibility is that the meeting point $m$ belongs to $\Upsilon$ or that it does not exist at all. (If it belongs to $\Upsilon$ it means that $\muu$ joins a $\g_j$.)

   \paragraph{Proofs of Corollaries \ref{convex-concave} and \ref{median}.} Consider an arc of a streamline $\be$   from the boundary point $ x \in \partial \Omega$ to the first meeting point $y$ on some $\g_j$. Now the speed is constant, say $|\nabla u(\be(t))| =c,\, 0 \leq t \leq T$ and so
   $$S\,=\int_0^T\!|\nabla u(\be(t))|\,dt\,=\,cT,\quad u(y) \,=\int_0^T\!|\nabla u(\be(t))|^2\,dt\,=\,c^2T\,=\,cS$$
   and the statement about the length follows. To prove the convexity of this arc, we suppose that the initial point lies  on the segment $P_kM_k$. If the arc is not convex, we can find a chord   with endpoints $a$ and $b$ on $\be$ so that the chord lies between $\be$ and the streamline from the corner $P_k$. But now $|\nabla u| \leq c$ on the chord $\overline{ab}$. This is because every streamline intersecting the chord starts somewhere at the boundary segment $P_k\,x$ and so its initial speed is smaller than the speed  $|\nabla u(x)|\,=\,c$. (The monotonicity of $|\nabla u|$ along the side  is needed.) Recall that the speed is constant along the streamline  $\be$. Thus
   
   \begin{gather*}
     u(b)-u(a) = \int_{t_a}^{t_b}\!|\nabla u(\be(t))|^2\,dt = c\int_{t_a}^{t_b}\!|\nabla u(\be(t))|\,dt = cS_{a,b}\\
       u(b)-u(a) = \int_0^1\!\langle\nabla u(a+t(b-a)),b-a\rangle dt \leq c\,|b-a|,
   \end{gather*}
   but the arclength $S_{a,b} > |b-a|$ if the arc is not a straight line segment. This proves  Corollary \ref{convex-concave}.

   Now Corollary \ref{median} follows immediately, because a chord of the 
median can lie neither  to the left nor to the right of it, since the median  has higher speed than the streamlines on both sides of it. \qquad $\Box$

\bigskip
{\small \paragraph{Acknowledgements.}  We thank Juan Manfredi at the University of Pittsburgh for a careful reading of the original manuscript. Erik Lindgren was supported by the Swedish Research Council, grant no. 2017-03736. Peter Lindqvist was supported by The Norwegian Research Council, grant no. 250070 (WaNP).}
   
\bigskip

\noindent 
{\small{\textsf{Erik Lindgren\\  Department of Mathematics\\ Uppsala University\\ Box 480\\
751 06 Uppsala, Sweden}  \\
\textsf{e-mail}: erik.lindgren@math.uu.se\\}

\bigskip
\noindent {\small\textsf{Peter Lindqvist\\ Department of
   Mathematical Sciences\\ Norwegian University of Science and
  Technology\\ N--7491, Trondheim, Norway}\\
\textsf{e-mail}: peter.lindqvist@ntnu.no}

\end{document}